\documentclass{amsart}
\usepackage{eurosym}
\usepackage{amssymb}
\usepackage{amsmath}
\usepackage{amstext}
\usepackage{amsfonts}
\usepackage{amsthm}
\usepackage{mathrsfs}
\usepackage{centernot}
\usepackage{enumitem}
\usepackage{graphicx}

\newtheorem*{theorem*}{Theorem}

\newtheorem{thmx}{Theorem}

\newtheorem{lem}[thmx]{Lemma}

\newtheorem{corollary}{Corollary}

\newtheorem{remark}{Remark}

\begin{document}
\title[PGT on modular surface]{Prime geodesic theorem for the modular surface%
}
\author{Muharem Avdispahi\'{c}}
\address{University of Sarajevo, Department of Mathematics, Zmaja od Bosne
33-35, 71000 Sarajevo, Bosnia and Herzegovina}
\email{mavdispa@pmf.unsa.ba}

\begin{abstract}
Under the generalized Lindel\"{o}f hypothesis, the exponent in the error
term of the prime geodesic theorem for the modular surface is reduced to $%
\frac{5}{8}+\varepsilon$ outside a set of finite logarithmic measure.
\end{abstract}

\maketitle

\section{Introduction}

Let $\Gamma =PSL\left( 2,%
\mathbb{Z}
\right) $ be the modular group and $\mathcal{H}$ the upper half-plane
equipped with the hyperbolic metric. The norms $N(P_{0})$ of primitive
conjugacy classes $P_{0}$ in $\Gamma $ are sometimes called pseudo-primes.
The length of the primitive closed geodesic on the modular surface $\Gamma
\setminus \mathcal{H}$ joining two fixed points, which are the same for all
representatives of $P_{0}$, equals $\log (N(P_{0}))$. The statement about
the number $\pi _{\Gamma }(x)$ of classes $P_{0}$ such that $N(P_{0})\leq x$%
, for $x>0$, is known as the prime geodesic theorem, PGT.

The main tool in the proof of PGT is the Selberg zeta function, defined by
\begin{equation*}
Z_{\Gamma }(s)=\underset{\{P_{0}\}}{\prod }\underset{k=0}{\overset{\infty }{%
\prod }}(1-N(P_{0})^{-s-k})\text{, }\text{Re}(s)>1\text{,}
\end{equation*}
and meromorphicaly continued to the whole complex plane.

The relationship between the prime geodesic theorem and the distribution of
zeros of the Selberg zeta function resembles to a large extent the
relationship between the prime number theorem and the zeros of the Riemann
zeta.

However, the function $Z_{\Gamma }$ satisfies the Riemann hypothesis. It is
an outstanding open problem whether the error term in the prime geodesic
theorem is $O(x^{\frac{1}{2}+\varepsilon})$ as it would be the case in the
prime number theorem once the Riemann hypothesis be proved.

The obstacles in establishing an analogue of von Koch's theorem \cite[p. 84]%
{I:1932} in this setting comes from the fact that $Z_{\Gamma}$ is a
meromorphic function of order $2$, while the Riemann zeta is of order $1$.

In the case of Fuchsian groups $\Gamma \subset PSL\left( 2,%
\mathbb{R}
\right) $, the best estimate of the remainder term in PGT is still $O\left(
\frac{x^{\frac{3}{4}}}{\log x}\right)$ obtained by Randol \cite{R:1977} (see
also \cite{B:1992}, \cite{A:2018} for different proofs). We note that its
analogue $O\left( x^{\frac{3}{2} d_{0}}\left( \log x\right) ^{-1}\right) $
is valid also for strictly hyperbolic manifolds of higher dimensions, where $%
d_{0}=\frac{d-1}{2}$ and $d\geq 3$ is the dimension of a manifold \cite[%
Theorem 1]{AG:2012}.

The attempts to reduce the exponent $\frac{3}{4}$ in PGT were successful
only in special cases. The chronological list of improvements for the
modular group $\Gamma = PSL(2,\mathbb{Z})$ includes $\frac{35}{48}%
+\varepsilon $ (Iwaniec \cite{I2:1984}), $\frac{7}{10}+\varepsilon $ (Luo
and Sarnak \cite{LS:1995}), $\frac{71}{102}+\varepsilon $ (Cai \cite{C:2002}%
) and the present $\frac{25}{36}+\varepsilon $ (Soundararajan and Young \cite%
{SY:2013}).

Iwaniec \cite{I1:1984} remarked that the generalized Lindel\"{o}f hypothesis
for Dirichlet $L$-functions would imply $\frac{2}{3}+\varepsilon $.

We proved \cite{A1:2017} that $\frac{2}{3}+\varepsilon$ is valid outside a
set of finite logarithmic measure. In the present note, we relate the error
term in the Gallagherian $PGT$ on $PSL(2,\mathbb{Z})$ to the subconvexity
bound for Dirichlet $L$- functions. This enables us to replace $\frac{2}{3}%
+\varepsilon$ by $\frac{5}{8}+\varepsilon$ under the generalized Lindel\"{o}%
f hypothesis. More precisely, the main result of this paper is the following
theorem.

\begin{theorem*}
\label{thm1} Let $\Gamma = PSL(2,\mathbb{Z})$ be the modular group, $%
\varepsilon >0$ arbitrarily small and $\theta$ be such that
\begin{equation*}
L\left(\frac{1}{2}+it, \chi_{D}\right) \ll \left(1+\left\vert t \right\vert
\right)^{A} \left\vert D \right\vert ^{\theta +\varepsilon}
\end{equation*}%
for some fixed $A > 0$, where $D$ is a fundamental discriminant. There
exists a set $B$ of finite logarithmic measure such that
\begin{equation*}
\pi _{\Gamma }\left( x\right) = \int_{0}^{x}\frac{dt}{\log t}+O\left( x^{%
\frac{5}{8}+\frac{\theta}{4}+\varepsilon }\right)\ \ \left( x\rightarrow
\infty , x\notin B\right) \text{.}
\end{equation*}
\end{theorem*}

Inserting the Conrey-Iwaniec \cite{C1:2000} value $\theta =\frac{1}{6}$ into
Theorem, we obtain

\begin{corollary}
\label{cor1}
\begin{equation*}
\pi _{\Gamma }\left( x\right) = \mathit{li} \left( x\right) +O\left( x^{%
\frac{2}{3}+\varepsilon }\right)\ \ \left( x\rightarrow \infty , x\notin
B\right) \text{.}
\end{equation*}
\end{corollary}

Any improvement of $\theta$ immediately results in the obvious improvement
of the error term in PGT. Taking into account that the Lindel\"{o}f
hypothesis allows $\theta=0$, we get

\begin{corollary}
\label{cor2} Under the Lindel\"{o}f hypothesis,
\begin{equation*}
\pi _{\Gamma }\left( x\right) = \mathit{li} \left( x\right) +O\left( x^{%
\frac{5}{8}+\varepsilon }\right)\ \ \left( x\rightarrow \infty , x\notin
B\right) \text{.}
\end{equation*}
\end{corollary}

\begin{remark}
The obtained exponent for strictly hyperbolic Fuchsian groups is $\frac{7}{10%
}+\varepsilon$ outside a set of finite logarithmic measure \cite{A2:2017}
and coincides with the above mentioned Luo-Sarnak unconditional result for $%
\Gamma = PSL(2,\mathbb{Z})$. In the case of a cocompact Kleinian group or a
noncompact congruence group for some imaginary quadratic number field, the
respective Gallagherian bound is $\frac{13}{9}+\varepsilon$ \cite{A3:2017}.
\end{remark}

\section{Preliminaries.}

The motivation for Theorem comes from several sources, including Gallagher
\cite{G1:1980}, Iwaniec \cite{I2:1984} and Balkanova and Frolenkov \cite%
{BF:2018}.

Recall that $\pi _{\Gamma }\left( x\right) =\mathit{li}\left( x\right)
+O\left( x^{\frac{5}{8}+\frac{\theta }{4}+\varepsilon }\right) $ is
equivalent to $\psi _{\Gamma }\left( x\right) =x+O\left( x^{\frac{5}{8}+%
\frac{\theta }{4}+\varepsilon }\right) $, where $\psi _{\Gamma }\left(
x\right) =\underset{N\left( P_{0}\right) ^{k}\leq x}{\sum }\log N\left(
P_{0}\right) $ is the $\Gamma $ analogue of the classical Chebyshev function
$\psi $.

Under the Riemann hypothesis, Gallagher improved von Koch's remainder term
in the prime number theorem from $\psi(x)=x+O\left(x^{\frac{1}{2}}(\log
x)^{2}\right)$ to $\psi(x)=x+O\left(x^{\frac{1}{2}}(\log \log x)^{2}\right)$
outside a set of finite logarithmic measure.

Following Koyama \cite{K:2016}, we shall apply the next lemma \cite{G:1970}
due to Gallagher to our setting.

\begin{lem}
\label{lemma1} Let $A$ be a discrete subset of $\mathbb{R}$ and $\eta \in
(0,1)$. For any sequence $c(\nu )\in \mathbb{C}$, $\nu \in A$, let the
series
\begin{equation*}
S\left( u\right) =\underset{\nu \in A}{\sum }c\left( \nu \right) e^{2\pi
i\nu u}
\end{equation*}%
be absolutely convergent. Then%
\begin{equation*}
\int_{-U}^{U}\left\vert S\left( u\right) \right\vert ^{2}du\leq \left( \frac{%
\pi \eta }{\sin \pi \eta }\right) ^{2}\int_{-\infty }^{+\infty }\left\vert
\frac{U}{\eta }\underset{t\leq \nu \leq t+\frac{\eta }{U}}{\sum }c\left( \nu
\right) \right\vert ^{2}dt\text{.}
\end{equation*}
\end{lem}

Iwaniec \cite{I2:1984} established the following explicit formula with an
error term for $\psi _{\Gamma }$ on $\Gamma = PSL(2,\mathbb{Z})$.

\begin{lem}
\label{lemma2} For $1\leq T\leq \frac{x^{\frac{1}{2}}}{\left(\log x\right)^2}
$, one has
\begin{equation*}
\psi_{\Gamma} \left( x\right) =x+\underset{\left\vert \gamma \right\vert
\leq T}{\sum }\frac{x^{\rho }}{\rho } +O\left(\frac{x}{T} \left(\log
x\right)^2\right)\text{,}
\end{equation*}
where $\rho = \frac{1}{2}+i\gamma$ denote zeros of $Z_{\Gamma}$.
\end{lem}

Recently, O. Balkanova and D. Frolenkov have proved the following estimate.

\begin{lem}
\label{lemma3}
\begin{eqnarray*}
\underset{\left\vert \gamma \right\vert \leq Y}{\sum }x^{i\gamma } &\ll
&\max \left( x^{\frac{1}{4}+\frac{\theta }{2}}Y^{\frac{1}{2}},x^{\frac{%
\theta }{2}}Y\right) \log ^{3}Y\text{,} \\
\underset{\left\vert \gamma \right\vert \leq Y}{\sum }x^{i\gamma } &\ll
&Y\log ^{2}Y\text{ if }Y>\frac{x^{\frac{1}{2}+\frac{7}{6}\theta }}{\kappa
\left( x\right) }\text{,}
\end{eqnarray*}%
where $\rho =\frac{1}{2}+i\gamma $ are the zeros of $Z_{\Gamma }$, $\theta $
is the subconvexity exponent for Dirichlet $L-$functions, and $\kappa \left(
x\right) $ is the distance from $\sqrt{x}+\frac{1}{\sqrt{x}}$ to the nearest
integer.
\end{lem}

\section{Proof of Theorem.}

Inserting $T=\frac{x^{\frac{1}{2}}}{\left( \log x\right) ^{2}}$ into Lemma %
\ref{lemma2}, we obtain
\begin{equation}
\psi _{\Gamma }\left( x\right) =x+\underset{\left\vert \gamma \right\vert
\leq T}{\sum }\frac{x^{\rho }}{\rho }+O\left( x^{\frac{1}{2}}\left( \log
x\right) ^{4}\right) \text{.}  \label{eq1}
\end{equation}

We would like to bound the expression $\underset{\left\vert \gamma
\right\vert \leq T}{\sum }\frac{x^{\rho }}{\rho }$, where $Y\in \left(
0,T\right) $ is a parameter to be determined later on.

Let $n=\left\lfloor \log x\right\rfloor $ and $B_{n}=\left\{ x\in \left[
e^{n},e^{n+1}\right) :\left\vert \underset{\left\vert \gamma \right\vert
\leq T}{{\sum }}\frac{x^{i\gamma }}{\rho }\right\vert >x^{\varepsilon }Y^{%
\frac{1}{2}}\right\} $. Looking at the logarithmic measure of $B_{n}$, we get%
\begin{eqnarray}
\mu ^{\ast }B_{n} &=&\underset{B_{n}}{\int }\frac{dx}{x}=\underset{A_{n}}{%
\int }x^{2\varepsilon }Y\frac{dx}{x^{1+2\varepsilon }Y}\leq \overset{e^{n+1}}%
{\underset{e^{n}}{\int }}\left\vert \underset{\left\vert \gamma \right\vert
\leq Y}{{\sum }}\frac{x^{i\gamma }}{\rho }\right\vert ^{2}\frac{dx}{%
x^{1+2\varepsilon }Y}  \label{eq2} \\
&\leq &\frac{1}{e^{2n\varepsilon }Y}\overset{e^{n+1}}{\underset{e^{n}}{\int }%
}\left\vert \underset{\left\vert \gamma \right\vert \leq Y}{{\sum }}\frac{%
x^{i\gamma }}{\rho }\right\vert ^{2}\frac{dx}{x}\text{.}  \notag
\end{eqnarray}

After substitution $x=e^{n}\cdot e^{2\pi \left( u+\frac{1}{4\pi }\right) }$,
the last integral becomes
\begin{equation*}
2\pi \overset{\frac{1}{4\pi }}{\underset{-\frac{1}{4\pi }}{\int }}\left\vert
\underset{\left\vert \gamma \right\vert \leq T}{{\sum }}\frac{e^{\left( n+%
\frac{1}{2}\right) i\gamma }}{\rho }e^{2\pi i\gamma u}\right\vert ^{2}du%
\text{.}
\end{equation*}

Applying Lemma \ref{lemma1}, with $\eta =U=\frac{1}{4\pi }$ and $c_{\gamma }=%
\frac{e^{\left( n+\frac{1}{2}\right) i\gamma }}{\rho }$ for $\left\vert
\gamma \right\vert \leq T$, $c_{\gamma }=0$ otherwise, we get
\begin{equation}
\overset{\frac{1}{4\pi }}{\underset{-\frac{1}{4\pi }}{\int }}\left\vert
\underset{\left\vert \gamma \right\vert \leq T}{{\sum }}\frac{e^{\left( n+%
\frac{1}{2}\right) i\gamma }}{\rho }e^{2\pi i\gamma u}\right\vert ^{2}du\leq
\left( \frac{\frac{1}{4}}{\sin \frac{1}{4}}\right) ^{2}\overset{+\infty }{%
\underset{-\infty }{\int }}\left( \sum_{\substack{t<\gamma \leq t+1 \\
\left\vert \gamma \right\vert \leq Y}}^{+\infty} \frac{1}{\left\vert \rho \right\vert }\right) ^{2}dt\text{.}
\label{eq3}
\end{equation}

Note that $\sum_{\substack{t<\gamma \leq t+1}}\frac{1}{\left\vert \rho \right\vert }=O\left( 1\right) $ since $\#\left\{ \gamma
:t<\left\vert \gamma \right\vert \leq t+1\right\} =O\left( t\right) $ by the
Weyl law.

Thus,%
\begin{equation}
\overset{+\infty }{\underset{-\infty }{\int }}\left( \sum_{\substack{t<\gamma \leq t+1 \\
\left\vert \gamma \right\vert \leq Y}}^{+\infty} \frac{1}{\left\vert \rho \right\vert }\right) ^{2}dt=O\left(\int\limits_{0}^{Y}dt\right) =O\left( Y\right) \text{.}  \label{eq4}
\end{equation}

The relations (\ref{eq2}), (\ref{eq3}) and (\ref{eq4}) imply $\mu ^{\ast
}B_{n}\ll \frac{Y}{e^{2n\varepsilon }Y}=\frac{1}{e^{2n\varepsilon }}$.
Hence, the set $B=\cup B_{n}$ has a finite logarithmic measure.

For $x\notin B$, we have $\left\vert \underset{\left\vert \gamma \right\vert
\leq Y}{{\sum }}\frac{x^{i\gamma }}{\rho }\right\vert \leq x^{\varepsilon
}Y^{\frac{1}{2}}$, i.e.
\begin{equation}
\left\vert \underset{\left\vert \gamma \right\vert \leq Y}{{\sum }}\frac{%
x^{\rho }}{\rho }\right\vert \leq x^{\frac{1}{2}+\varepsilon }Y^{\frac{1}{2}}%
\text{.}  \label{eq5}
\end{equation}

Now, we rely on Lemma \ref{lemma3} to estimate $\left\vert \underset{%
Y<\left\vert \gamma \right\vert \leq T}{{\sum }}\frac{x^{\rho }}{\rho }%
\right\vert $. Let us put $S\left( x,T\right) =\underset{\left\vert \gamma
\right\vert \leq T}{{\sum }}x^{i\gamma }$. By Abel's partial summation, we
have%
\begin{equation*}
\underset{Y<\left\vert \gamma \right\vert \leq T}{{\sum }}\frac{x^{i\gamma }%
}{\rho }=\frac{S\left( x,T\right) }{\frac{1}{2}+iT}-\frac{S\left( x,Y\right)
}{\frac{1}{2}+iY}+i\int\limits_{Y}^{T}\frac{S\left( x,u\right) }{\left(
\frac{1}{2}+iu\right) ^{2}}du\text{.}
\end{equation*}

Multiplying the last relation by $x^{\frac{1}{2}}$ and recalling that Lemma %
\ref{lemma3} yields $\underset{\left\vert \gamma \right\vert \leq Y}{\sum }%
x^{i\gamma }\ll x^{\frac{1}{4}+\frac{\theta }{2}+\varepsilon }Y^{\frac{1}{2}}$ for $Y<T=%
\frac{x^{\frac{1}{2}}}{\left( \log x\right) ^{2}}$, we get

\begin{equation}
\left\vert \underset{Y<\left\vert \gamma \right\vert \leq T}{{\sum }}\frac{%
x^{\rho }}{\rho }\right\vert \ll \frac{x^{\frac{3}{4}+\frac{\theta }{2}%
+\varepsilon }}{T^{\frac{1}{2}}}+\frac{x^{\frac{3}{4}+\frac{\theta }{2}%
+\varepsilon }}{Y^{\frac{1}{2}}}+\int\limits_{Y}^{T}\frac{x^{\frac{3}{4}+%
\frac{\theta }{2}+\varepsilon }u^{\frac{1}{2}}}{u^{2}}du\ll \frac{x^{\frac{3%
}{4}+\frac{\theta }{2}+\varepsilon }}{Y^{\frac{1}{2}}}\text{.}  \label{eq6}
\end{equation}

Combining (\ref{eq5}) and (\ref{eq6}), we see that the optimal choice for the
parameter $Y$ is $Y\approx x^{\frac{1}{4}+\frac{\theta }{2}}$. Then, $%
\underset{\left\vert \gamma \right\vert \leq T}{{\sum }}\frac{x^{\rho }}{%
\rho }=O\left( x^{\frac{1}{2}+\varepsilon }Y^{\frac{1}{2}}\right) =O\left(
x^{\frac{5}{8}+\frac{\theta }{4}+\varepsilon }\right) $ for $x\notin B$.

The relation (\ref{eq1}) becomes%
\begin{equation*}
\psi _{\Gamma }\left( x\right) =x+O\left( x^{\frac{5}{8}+\frac{\theta }{4}%
+\varepsilon }\right) \text{ }\left( x\rightarrow \infty \text{, }x\notin
B\right) \text{,}
\end{equation*}%
as asserted.


\bigskip

\begin{thebibliography}{99}
\bibitem{A:2018} Avdispahi\'{c}, M. ``On Koyama's refinement of the prime
geodesic theorem.'' \textit{Proc. Japan Acad. Ser. A} 94, no. 3 (2018),
21--24.

\bibitem{A1:2017} Avdispahi\'{c}, M. ``Gallagherian $PGT$ on $PSL(2,\mathbb{Z%
})$.'' \textit{Funct. Approximatio. Comment. Math.} doi:10.7169/facm/1686

\bibitem{A2:2017} Avdispahi\'{c}, M. ``Prime geodesic theorem of Gallagher
type.'' arXiv:1701.02115.

\bibitem{A3:2017} Avdispahi\'{c}, M. ``On the prime geodesic theorem for
hyperbolic $3$-manifolds.'' \textit{Math. Nachr.} (to appear; cf.
arXiv:1705.05626).

\bibitem{AG:2012} Avdispahi\'{c}, M., and D\v{z}. Gu\v{s}i\'{c}. ``On the
error term in the prime geodesic theorem.'' \textit{Bull. Korean Math. Soc.}
49, no. 2 (2012), 367--372.

\bibitem{BF:2018} Balkanova, O., and D. Frolenkov. ``Bounds for the spectral
exponential sum.'' arXiv:1803.04201.

\bibitem{B:1992} Buser, P. \textit{Geometry and spectra of compact Riemann
surfaces}, Progress in Mathematics, Vol. 106, Birkh\"{a}user,
Boston-Basel-Berlin, 1992.

\bibitem{C:2002} Cai, Y. ``Prime geodesic theorem.'' \textit{J. Th\'{e}or.
Nombres Bordeaux} 14, no. 1 (2002), 59--72.

\bibitem{C1:2000} Conrey, J. B. and H. Iwaniec. ``The cubic moment of
central values of automorphic L-functions.'' \textit{Ann. of Math. (2)} 151,
no. 3 (2000), 1175--1216.

\bibitem{G:1970} Gallagher, P. X. ``A large sieve density estimate near $%
\sigma =1$.'' \textit{Invent. Math.} 11 (1970), 329--339.

\bibitem{G1:1980} Gallagher, P. X. ``Some consequences of the Riemann
hypothesis.'' \textit{Acta Arith.} 37 (1980), 339--343.

\bibitem{H:1976} Hejhal, D. A. \textit{The Selberg trace formula for $%
PSL(2,R)$. Vol I}, Lecture Notes in Mathematics, Vol 548, Springer, Berlin,
1976.

\bibitem{I:1932} Ingham, A. E. \textit{The distribution of prime numbers},
Cambridge University Press, 1932.

\bibitem{I1:1984} Iwaniec, H. ``Non-holomorphic modular forms and their
applications.'' In \textit{Modular forms} (Durham, 1983), 157-\^{a}\euro
``196, Ellis Horwood Ser. Math. Appl.: Statist. Oper. Res., Horwood,
Chichester, 1984.

\bibitem{I2:1984} Iwaniec, H. ``Prime geodesic theorem.'' \textit{J. Reine
Angew. Math.} 349 (1984), 136--159.

\bibitem{K:2016} Koyama, S. ``Refinement of prime geodesic theorem.''
\textit{Proc. Japan Acad. Ser A Math. Sci.} 92, no. 7 (2016), 77--81.

\bibitem{LS:1995} Luo, W. and P. Sarnak. ``Quantum ergodicity of
eigenfunctions on $PSL_{2}(Z)\backslash H^{2}$.'' \textit{Inst. Hautes \'{E}%
tudes Sci. Publ. Math.} no. 81 (1995), 207--237.

\bibitem{R:1977} Randol, B. ``On the asymptotic distribution of closed
geodesics on compact Riemann surfaces.'' \textit{Trans. Amer. Math. Soc.}
233 (1977), 241--247.

\bibitem{SY:2013} Soundararajan, K. and M. P. Young. ``The prime geodesic
theorem.'' \textit{J. Reine Angew. Math.} 676 (2013), 105--120.
\end{thebibliography}
\end{document}